\documentclass{article}

\usepackage[english]{babel}

\usepackage[letterpaper,top=2cm,bottom=2cm,left=3cm,right=3cm,marginparwidth=1.75cm]{geometry}

\usepackage{verbatim}
\usepackage{caption}
\usepackage{float}
\usepackage{pifont}
\usepackage{enumerate}
\usepackage{amsmath}
\usepackage{graphicx}
\usepackage[colorlinks=true, allcolors=blue]{hyperref}
\newtheorem{theorem}{Theorem}[section]
\newtheorem{conjecture}[theorem]{Conjecture}
\newtheorem{lemma}[theorem]{Lemma}
\newtheorem{corollary}[theorem]{Corollary}

\def\qed{\hfill \rule{4pt}{7pt}}

\def\pf{\noindent {\it Proof. }} 

\title{Graphs with girth 9 and without longer odd holes are 3-colorable}
\author{Yan Wang$^{1,}$ \thanks{Supported by National Key R\&D Program of China under Grant No. 2022YFA1006400 and Explore-X project of SJTU, email: yan.w@sjtu.edu.cn},  \;\; Rong Wu$^{2,}$ \thanks{Supported by National Key R\&D Program of China under Grant No. 2022YFA1006400, email: wurong@sjtu.edu.cn} \\\\
    \small $^1$School of Mathematical Sciences, CMA-Shanghai\\
    \small Shanghai Jiao Tong University, 800 Dongchuan Road, Shanghai 200240, China\\
    \small $^2$School of Mathematical Sciences \\
    \small Shanghai Jiao Tong University, 800 Dongchuan Road, Shanghai 200240, China}

\begin{document}
\maketitle

\begin{abstract}
For a number $l\geq 2$, let ${\cal{G}}_l$ denote the family of graphs which have girth $2l+1$ and have no odd hole with length greater than $2l+1$. Wu, Xu and Xu conjectured that every graph in $\bigcup_{l\geq 2} {\cal{G}}_{l}$ is $3$-colorable. Chudnovsky et al., Wu et al., and Chen showed that every graph in ${\cal{G}}_2$, ${\cal{G}}_3$ and $\bigcup_{l\geq 5} {\cal{G}}_{l}$ is $3$-colorable respectively. 
In this paper, we prove that every graph in ${\cal{G}}_4$ is $3$-colorable.
This confirms Wu, Xu and Xu's conjecture.
\end{abstract}

\section{Introduction}

All graphs considered in this paper are finite, simple, and undirected. 
Let $G$ be a graph and let $S$ be a subset of $V(G)$.
We use $G[S]$ to denote the subgraph of $G$ induced by $S$.
Let $x\in V(G)$ (we also write $x\in G$ if there is no confusion), we use $N_S(x)$ to denote the neighbours of $x$ in $S$.
For subgraphs $H$ and $H'$ of $G$, we use $H \triangle H'$ to denote the symmetric difference of $H$ and $H'$, that is, 
$V(H\triangle H')=V(H)\cup V(H')\backslash \{V(H)\cap V(H')\}$ and $E(H\triangle H')=E(H)\cup E(H')\backslash \{E(H)\cap E(H')\}$.

Let $P$ be an $(x, y)$-path, that is, the ends of $P$ are $x$ and $y$.
And we usually use $xPy$ to denote $P$, if there is no confusion, we just write it as $P$.
Let $P^*$ denote the set of internal vertices of $P$.
Let $C$ be a cycle and $u, v$ be two vertices of $C$. We use $C(u, v)$ to denote the subpath of $C$ from $u$ to $v$ in clockwise order, and $C^*(u, v)$ to denote the set of internal vertices of $C(u, v)$.

A graph $G$ is $k$-colorable if there exists a mapping $c: V(G)\rightarrow \{1, 2, \cdots, k\}$ such that $c(u)\neq c(v)$ whenever $uv\in E(G)$.
The \textit{chromatic number $\chi(G)$} of $G$ is the minimum integer $k$ such that $G$ is $k$-colorable.
The clique number $\omega(G)$ of $G$ is the maximum integer $k$ such that $G$ contains a complete graph of size $k$.
For a graph $G$, if $\chi(G)=\omega(G)$, then we call $G$ a \textit{perfect} graph.
For a graph $H$, we say that $G$ is \textit{$H$-free} if $G$ induces no $H$ (i.e., $G$ has no induced subgraph isomorphic to $H$).
Let $\cal F$ be a family of graphs. We say that $G$ is \textit{$\cal F$-free} if $G$ induces no member of $\cal F$.
If there exists a function $\phi$ such that $\chi(G)\leq \phi(\omega(G))$ for each $G\in {\cal F}$, then we say that $\cal{F}$ is $\chi$-\textit{bounded class}, and call $\phi$ a \textit{binding function} of $\cal{F}$.
The concept of $\chi$-boundedness was put forward by Gy\'{a}rf\'{a}s in 1975 \cite{g2}.
Studying which family of graphs can be $\chi$-bounded, and finding the optimal binding function for $\chi$-bounded class are important problems in this area.
Since clique number is a trivial lower bound of chromatic number, if a family of $\chi$-bounded graphs has a linear binding function, then the linear function must be asymptotically optimal binding function of this family.
For more recent information on $\chi$-bounded problems, see \cite{ss}.

A hole in a graph is an induced cycle of length at least $4$.
A hole is said to be \textit{odd} (resp. \textit{even}) if it has odd (resp. even) length.
Addario-Berry, Chudnovsky, Havet, Reed and Seymour \cite{achrs}, and Chudnovsky and Seymour \cite{cs08}, proved that every even hole free graph has a vertex whose neighbours are the union of two cliques, which implies that $\chi(G)\leq 2\omega(G)-1$.
However, the situation becomes much more complicated on odd hole free graphs.
The Strong Perfect Graph Theorem \cite{cs02} asserts that a graph is perfect if and only if it induces neither odd holes nor their complements.
Confirming a conjecture of Gy\'{a}rf\'{a}s \cite{g2},
Scott and Seymour \cite{sso} proved that odd hole free graphs are $\chi$-bounded with binding function $\frac{2^{2^{\omega(G)+2}}}{48(\omega(G)+2)}$.
Ho\'{a}ng and McDiarmid \cite{hm} conjectured for an odd hole free graph $G$, $\chi(G)\leq 2^{\omega(G)-1}$.
A graph is said to be \textit{short-holed} if every hole of it has length 4. Sivaraman \cite{sivaraman} conjectured that $\chi(G)\leq \omega^2(G)$ for all short-holed graphs whereas the best known upper-bound is $\chi(G)\leq 10^{20}2^{\omega^2(G)}$ due to Scott and Seymour \cite{ss}.

The girth of a graph $G$, denoted by $g(G)$, is the minimum length of a cycle in $G$.
Let $l\geq 2$ be an integer. Let ${\cal G}_l$ denote the family of graphs that have girth $2l+1$ and have no odd holes of length at least $2l+3$. 
The graphs in ${\cal G}_2$ are called \textit{pentagraphs}, and the graphs in ${\cal G}_3$ are called \textit{heptagraphs}.
Robertson \cite{NPRZ11} conjectured  that the Petersen graph is the only non-bipartite pentagraph which is $3$-connected and internally $4$-connected. Plummer and Zha \cite{PZ14} presented some counterexamples to Robertson's conjecture, and conjectured that every pentagraph is $3$-colorable.
Xu, Yu and Zha \cite{xyz17} proved that every pentagraph is $4$-colorable.
Generalizing the result of \cite{xyz17}, Wu, Xu and Xu \cite{wxx22} proved that graphs in $\bigcup_{l\geq 2}{\cal G}_l$ are $4$-colorable and proposed the following conjecture.
\begin{conjecture} \label{wxx22}{\em{\cite{wxx22}}}
Graphs in $\bigcup_{l\geq 2}{\cal G}_l$ are $3$-colorable.
\end{conjecture}

Recently, Chudnovsky and Seymour \cite{cs22} confirmed that pentagraphs are $3$-colorable.
Wu, Xu and Xu \cite{wxx22+} showed that heptagraphs are $3$-colorable.
More recently, Chen \cite{c23} proved that all graphs in $\bigcup_{l\geq 5}{\cal G}_l$ are $3$-colorable.
In this paper, we prove Conjecture \ref{wxx22}.

\begin{theorem} \label{main theorem}
Graphs in ${\cal G}_4$ are $3$-colorable.
\end{theorem}

\section{Preliminary}

In this section, we collect some useful lemmas. 
The authors of \cite{cs22} proved the following lemma for $l=2$, but their proof also works for $l \ge 2$.

\begin{lemma}\label{cs} {\em \cite{cs22}}
For any number $l\geq 2$, every $4$-vertex-critical graph in ${\cal G}_l$ has neither $K_2$-cut or $P_3$-cut. 
\end{lemma}

\begin{lemma} \label{cz2.2}{\em \cite{cz22}}
For any number $k\geq 4$, each $k$-vertex-critical graph has no $2$-edge-cut. 
\end{lemma}

A theta graph is a graph that consists of a pair of distinct vertices joined by three internally disjoint paths.
Let $C$ be a hole of a graph $G$.
A path $P$ of $G$ is a chordal path of $C$ if $C\cup P$ is an induced theta-subgraph of $G$.

\begin{lemma} \label{cz2.3}{\em \cite{cz22}}
Let $l\geq 4$ be an integer and $C$ be an odd hole of a graph $G\in {\cal G}_l$.
Let $P$ be a chordal path of $C$, and $P_1$, $P_2$ be the internally disjoint paths of $C$ that have the same ends as $P$.
Assume that $|P|$ and $|P_1|$ have the same parity. Then
\begin{equation*}
 \mbox{$|P_1|=1$ or $l\geq |P_2|<|P_1|=|P|\geq l+1.$}  
\end{equation*}
In particular, when $|P_1|\geq 2$, all chordal paths of $C$ with the same ends as $P_1$ have length $|P_1|$.
\end{lemma}

\begin{lemma} \label{cz2.4}{\em \cite{cz22}}
Let $l \ge 4$ be an integer and  $x,y$ be the vertices of a graph $G\in{\cal G}_l$.
Let $X$ be a vertex cut of $G$ with $\{x, y\}\subseteq X\subseteq N[\{x, y\}]$, and $G_1$ be an induced subgraph of $G$ whose vertex set consists of $X$ and the vertex set of a component of $G-X$.
If all induced $(x, y)$-paths in $G_1$ have length $k$ with $4\leq k\leq l$, then $G$ has  a degree-2 vertex or a $K_2$-cut.
\end{lemma}

\begin{figure}[H]
\centering
\includegraphics[width=0.25\linewidth]{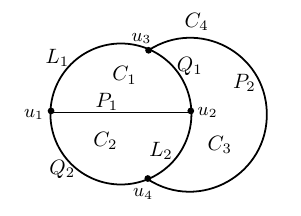}
\caption{\label{fig:$H$}$H$ is a $K_4$-subdivision. $u_1, u_2, u_3, u_4$ are the degree-3 vertices of $H$. $C_1, C_2, C_3, C_4$ are the faces cycles of $H$ and $\{P_1, P_2\}$, $\{Q_1, Q_2\}$, $\{L_1, L_2\}$ are the pairs of vertex disjoint arrises of $H$.}
\end{figure}

Let $H=(u_1,u_2,u_3,u_4,P_1,P_2,Q_1,Q_2,L_1,L_2)$ be a $K_4$-subdivision such that
$u_1, u_2, u_3, u_4$ are degree-3 vertices of $H$ 
and $P_1$ is a $(u_1, u_2)$-path, $P_2$ is a $(u_3, u_4)$-path, $Q_1$ is a $(u_2, u_3)$-path, $Q_2$ is a $(u_1, u_4)$-path, $L_1$ is a $(u_1, u_3)$-path, and $L_2$ is a $(u_2, u_4)$-path (see Figure \ref{fig:$H$}).
We call $P_1, P_2$, $Q_1, Q_2$, $L_1, L_2$ \textit{arrises} of $H$.
Let 
$C_1:=P_1\cup Q_1\cup L_1$, $C_2:=P_1\cup Q_2\cup L_2$, $C_3=P_2\cup Q_1\cup L_2$ and $C_4:=C_1\triangle C_2\triangle C_3$ be four holes in $H$. 
We call that $H$ is an \textit{odd $K_4$-subdivision} if $C_1$, $C_2$, $C_3$ and $C_4$ are odd holes.
If all arrises of an odd $K_4$-subdivision have the same length, then we call it a \textit{regular} odd $K_4$-subdivision, otherwise an \textit{irregular} odd $K_4$-subdivision.
If $C_1$ and $C_2$ are odd holes, $C_3$ and $C_4$ are even holes, $|Q_1|=1$ and $|L_2|\geq 2$, then we call $H$ a \textit{balanced $K_4$-subdivision of type $(1, 2)$}.

\begin{lemma} {\label{cz2.7}}{\em \cite{cz22}}
Let $l \ge 2$ be an integer and $H$ be a subgraph of a graph $G\in {\cal G}_l$. If $H$ is isomorphic to an odd $K_4$-subdivision, then the following statements hold.

{\em(1)} Each pair of vertex disjoint arrises have the same length and their lengths are at most $l$.

{\em(2)} $H$ is an induced subgraph of $G$.

{\em(3)} When $l\geq 3$, no vertex in $V(G)-V(H)$ has two neighbours in $H$.
\end{lemma}

Chen in \cite{c23} proved the following lemmas holds for $l\geq 5$. In fact, they are also true for $l=4$. 

\begin{lemma}{\label{c231}} {\em \cite{c23}}
Let $l\geq 4$ be an integer and $G$ be a graph in ${\cal G}_l$. If $G$ is $4$-vertex-critical, then $G$ does not contain a balanced $K_4$-subdivision of type $(1, 2)$.
\end{lemma}

Let $C$ be an odd hole of a graph $G$ and $s, t\in V(C)$ nonadjacent.
Let $P$ be an induced $(s, t)$-path.
If $V(C)\cap V(P^*)=\emptyset$, we call $P$ a jump or an $(s, t)$-jump over $C$.
Let $Q_1, Q_2$ be the internally disjoint $(s, t)$-paths of $C$.
If some vertex in $V(Q^*_1)$ has a neighbour in $V(P^*)$ and no vertex in $V(Q^*_2)$ has a neighbour in $V(P^*)$, we say that $P$ is a local jump over $C$ across $Q^*_1$.
When there is no need to point out $Q^*_1$, we will also say that $P$ is a local jump over $C$.
In particular, when $|V(Q^*_1)|=1$, we say that $P$ is a local jump over $C$ across one vertex.
When no vertex in $V(Q^*_1\cup Q^*_2)$ has  a neighbour in $V(P^*)$, we say that $P$ is a short jump over $C$. 
Hence, short jumps over $C$ are chordal paths of $C$. 
But chordal paths over $C$ maybe not short jumps over $C$ as the ends of a chordal path maybe adjacent.
When $P$ is a short jump over $C$, if $PQ_1$ is an odd hole, we say that $PQ_1$ is a jump hole over $C$ and $P$ is a short jump over $C$ across $Q^*_1$.

\begin{lemma}\label{c236}  {\em \cite{c23}}
Let $l\geq 4$ be an integer and $C$ be an odd hole of a graph $G\in {\cal G}_l$, $s, t$ be two vertices of $C$. 
Let $P$ be an $(s, t)$-jump, $Q_1=C(s, t)$ and $Q_2=C(t, s)$.
If $P$ is a local or short jump over $C$ across $Q^*_1$, 
then $|P|, |Q_2|$ have the same parity, and thus, $PQ_2$ is an even hole and $PQ_1$ is an odd cycle.
\end{lemma}

\begin{lemma}\label{c232}  {\em \cite{c23}}
Let $l\geq 4$ be an integer and $C$ be an odd hole of a graph $G\in {\cal G}_l$. 
If a jump $P$ over $C$ is not local, then $G[V(C\cup P)]$ contains a short jump over $C$.
\end{lemma}

\begin{lemma}\label{c233}  {\em \cite{c23}}
Let $ l\geq 4$ be an integer and $C$ be an odd hole of a graph $G\in {\cal G}_l$. 
If $P$ is a local $(v_1, v_2)$-jump over $C$,  then $G[V(C\cup P)]$ contains a local jump over $C$ across one vertex or a short jump over $C$.
\end{lemma}

It is easy to derive the following from Lemmas \ref{c232} and \ref{c233}.

\begin{corollary} \label{coro-jump}
Let $l\geq 4$ be an integer and $C$ be an odd hole of a graph $G\in {\cal G}_l$. 
If there is a jump over $C$, then $G$ contains either a short jump or a local jump over $C$ across one vertex.
\end{corollary}

Let $C$ be an odd hole and $P_i$ be an $(u_i, v_i)$-jump over $C$ for each integer $1\leq i\leq 2$.
If $u_1,v_1,u_2,v_2$ are disjoint vertices and $u_1,u_2,v_1,v_2$ appear on $C$ in this order, then we say that $P_1, P_2$ are \textit{crossing}; otherwise, they are \textit{uncrossing}. 
The following lemmas concern crossing jumps in $G$.

\begin{lemma}\label{c234}  {\em \cite{c23}}
Let $l\geq 4$ be an integer and $C$ be an odd hole of a graph $G\in {\cal G}_l$, and $P_i$ be a short $(u_i, v_i)$-jump over $C$ for each integer $1\leq i\leq 2$.
If $P_1$, $P_2$ are crossing, then $G$ contains an odd $K_4$-subdivision or a balanced $K_4$-subdivision of type $(1, 2)$.
\end{lemma}

\begin{lemma}\label{c235}  {\em \cite{c23}}
Let $l\geq 4$ be an integer and $C$ be an odd hole of a graph $G\in {\cal G}_l$.
For any integer $1\leq i\leq 2$, let $P_i$ be a short or local $(u_i, v_i)$-jump over $C$ such that $\{u_1, v_1\}\neq \{u_2, v_2\}$ and $u_1, u_2, v_2, v_1$ appear on $C$ in clockwise order. Assume that $P_1$ is across $C^*(v_1, u_1)$,  $P_2$ is across $C^*(u_2, v_2)$ and if $P_i$ is local, then $P_i$ is across one vertex for any integer $1\leq i\leq 2$. Then the following hold.
\begin{enumerate}[1)]
    \item 
    When $P_1, P_2$ are short, $G$ has an odd $K_4$-subdivision.
    \item 
    At least one of $P_i$ {\em(}$i\in [2]${\em)} is local, then at most two vertices in $V(C(v_1, u_1)\cup C(u_2, v_2))$ are not in a jump hole over $C$.
\end{enumerate}
\end{lemma}

\section{Proof of Theorem~\ref{main theorem}}

Let $H_1, H_2$ be vertex disjoint induced subgraphs of a graph $G$.
An induced $(v_1, v_2)$-path $P$ is a direct connection linking $H_1$ and $H_2$ if $v_1$ is the only vertex in $V(P)$ having a neighbour in $H_1$ and $v_2$ is the only vertex in $V(P)$ having a neighbour in $H_2$.

\begin{lemma} \label{odd lem}
Let $l\geq 4$ be an integer.
For each graph $G$ in ${\cal G}_l$, suppose $G$ is 4-vertex-critical, either $G$ has no odd $K_4$-subdivision or $G$ has an odd $K_4$-subdivision $H = (u_1,u_2,u_3,u_4,P_1,P_2,Q_1,Q_2,L_1,L_2)$ such that every minimal direct connection $(v_1, v_2)$-path linking $H\backslash {P^*_2}$ and $P^*_2$ must have $N_H(v_1)=N_{H\backslash {P^*_2}}(v_1)=\{u_3\}$ or $\{u_4\}$ and $N_H(v_2)=N_{{P^*_2}}(v_2)=N_{P_2^*}(N_H(v_1))$.
\end{lemma}
\pf Suppose $G$ has an odd $K_4$-subdivision denoted by $H$.
By Lemma \ref{cz2.7} (2), $H$ is an induced subgraph of $G$. Let $H=(u_1,u_2,u_3,u_4,P_1,P_2,Q_1,Q_2,L_1,L_2)$ (see Figure \ref{fig:$H$}) and $C_1:=P_1\cup Q_1\cup L_1$, $C_2:=P_1\cup Q_2\cup L_2$, $C_3=P_2\cup Q_1\cup L_2$ and $C_4:=C_1\triangle C_2\triangle C_3$. Since $H$ is an odd $K_4$-subdivision, $C_1$, $C_2$, $C_3$ and $C_4$ are odd holes. By Lemma \ref{cz2.7} (1), 
\begin{equation}\label{eq1}\tag{3.1}
    |P_1|=|P_2|\leq l, |Q_1|=|Q_2|\leq l, |L_1|=|L_2|\leq l. 
\end{equation}
Without loss of generality we may assume that $P_1$, $P_2$ are longest arrises in $H$.
Let $e, f$ be the edges of $P_2$ incident with $u_3, u_4$, respectively. Since $G$ is 4-vertex-critical, $\{e, f\}$ is not an edge-cut of $G$ by Lemma \ref{cs}.
So there exists a minimal path $P$ linking ${P^*_2}$ and $H \backslash V({P^*_2})$ such that $P$ is a $(v_1, v_2)$-path, $N_{H \backslash {P^*_2}}(v_1)\neq \emptyset$ and $N_{{P^*_2}}(v_2)\neq \emptyset$. 
Let $P$ be such a path.
By Lemma \ref{cz2.7} (3), let $N_{H \backslash {P^*_2}}(v_1)=\{x\}$ and $N_{{P^*_2}}(v_2)=\{y\}$.
Let $P'=xv_1Pv_2y$.
Note that $H\cup P'$ is induced by minimality of $P$.
Also note that any minimal path $P$ satisfies the following claims.

\medskip

\noindent \textbf{Claim 3.1.1} $x\notin \{u_1, u_2\}$.

\pf Suppose it is false. 
By symmetry we may assume that $x=u_1$.
Set $C'_4=u_3L_1u_1P'yP_2u_3$.
Since $C_4$ is an odd hole, by symmetry we may assume that $C'_4$ is an even hole and $C_4\triangle C'_4$ is an odd hole. Since $u_1P'yP_2u_3$ is a chordal path of $C_1$, by (\ref{eq1}) and Lemma \ref{cz2.3}, we have $|L_1|=1$.  So $|P_1|$=$|Q_1|=l$ by (\ref{eq1}) again. Since $C_4\triangle C'_4$ is an odd hole, $|u_1P'yP_2 u_4|=l+1$ which implies $|P'|\leq l$. Moreover, since $|P_2|=l$ and $|L_1|=1$, we have $|C'_4|\leq 2l$, a contradiction as $g(G) = 2l+1$. So $x\neq u_1$. \qed

\medskip

We say that $|P_1|-\min \{|Q_1|, |L_1|\}$ is the {\textit{difference}} of $H$.
Without loss of generality we may assume that among all odd $K_4$-subdivisions,
$H$ is chosen with smallest difference.

\medskip

\noindent \textbf{Claim 3.1.2} $x\notin V(P^*_1)$.

\pf Suppose it is false.
Without loss of generality, we may assume that $|L_1|\geq |Q_1|$. Set $C'_2=u_4Q_2u_1P_1xP'y$ $P_2u_4$. Since $C_4$ is an odd hole, either $C'_2$ or $C_4\triangle C'_2$ is an odd hole. When  $C_4\triangle C'_2$ is an odd hole, since $C_1\cup C_3\cup P'$ is an odd $K_4$-subdivision, by Lemma \ref{cz2.7} (1), $|P'|=|Q_1|$, $|u_1P_1x|=|u_4P_2y|$, and $|u_2P_1x|=|u_3P_2y|$. Since $|P'|=|Q_1|=|Q_2|$, we have $C'_2$ is an even hole of length $2(|Q_2|+|u_1P_1x|)$, implying $|L_1|+|u_1P_1x|\geq |Q_2|+|u_1P_1x|\geq l+1$. Then $C_4\triangle C'_2\triangle C_1$ is an even cycle of length at most $2l$, which is a contradiction. So $C'_2$ is an odd hole.

Since $C_2\cup C'_2\cup C_3$ is an odd $K_4$-subdivision, it follows from Lemma \ref{cz2.7} (1) that 

\begin{equation}\label{eq2}\tag{3.2}
    |P'|=|L_2|, |u_1P_1x|=|u_3P_2y|, |u_2P_1x|=|u_4P_2y|.
\end{equation}

Since $C_2\triangle C'_2\triangle C_1$ is an odd cycle of length larger than $2l+1$, it is not an odd hole, so
\begin{equation}\label{eq3}\tag{3.3}
    1\in \{|Q_2|, |u_2P_1x|, |u_3P_2y|\}
\end{equation}

When $|u_2P_1x|=1$, by (\ref{eq1}), (\ref{eq2}) and Lemma \ref{cz2.3}, we have $|L_1|=|P'|=l$. 
Moreover, as $l \ge |P_1|\ge |L_1| = l$, we have $|P_1| = l$.
So the difference of $H$ is $l-1$. 
Since $|u_2P_1x|$=$|u_4P_2y|$=$|Q_1|$=1, the graph $G[V(C_1\cup C'_2\cup P_2)]$ is an odd $K_4$-subdivision with difference $l-2$, 
which is a contradiction to the choice of $H$. So $|u_2P_1x|\geq 2$.
Assume that $|Q_2|=1$. Then $|L_1|=|P_1|=l$ by (\ref{eq1}).
Since $G[V(C'_2\cup C_2\cup C_3)]$ is an odd $K_4$-subdivision whose difference is at most $l-2$ as $|u_2P_1x|\geq 2$, which is contradiction to the choice of $H$. So $|Q_2|\geq 2$. Then $yu_3\in E(H)$ by (\ref{eq3}), implying $xu_1\in E(H)$ by (\ref{eq2}).
Then $|C_4\triangle C'_2|=2+2|L_1|$ by (\ref{eq1}) and (\ref{eq2}), so $|L_1|=l$ by (\ref{eq1}) again. Since $|P_1|\geq |L_1|$, we have $|P_1|=l$ and $|Q_1|=1$ by (\ref{eq1}), which is a contradiction as $|Q_2|\geq 2$.\qed

\medskip

\noindent \textbf{Claim 3.1.3} If $x\in \{u_3, u_4\}$, then $xy\in \{e, f\}$.

\pf By symmetry we may assume that $x=u_3$. 
Assume to the contrary that $x, y$ are non-adjacent.
Set $C'_3=u_3P'yP_2u_3$. Since $P'$ is a chordal path of $C_3$, 
we have that $C'_3$ is an odd hole by Lemma \ref{cz2.3} and (\ref{eq1}).
Since $C'_3\triangle C_3$ is an even hole, $C'_3\triangle C_3 \triangle C_2$ is an odd cycle and $|Q_1|=|L_2|=1$ by (\ref{eq1}) and Lemma \ref{cz2.3} again. Then $|P_1|=2l-1>l$ by (\ref{eq1}), which is a contradiction. So $e=xy$.\qed 

\medskip

\noindent \textbf{Claim 3.1.4} If $x\in V(L^*_1)$, then $|Q_1|=1$, $|P_1|=|L_1|=l$, $|P'|=2l-1$ and $xu_3, yu_3\in E(H)$.

\pf Set $C'_4=xL_1u_1Q_2u_4P_2yP'x$. 
Assume that $C_4\triangle C'_4$ is an even hole.
Since $x\neq u_3$, $xP'yP_2u_3$ is a chordal path of $C_1$. Hence, $xu_3\in E(H)$ by (\ref{eq1}) and Lemma \ref{cz2.3}.
Since $x\neq u_1$, the path $u_3xP'y$ is a chordal path of $C_3$.
If $yu_3\notin E(H)$, then $|u_3xP'y|=|u_3P_2y|\geq l+1$ by (\ref{eq1}) and Lemma \ref{cz2.3}, which is a contradiction to the fact that $|P_2|\leq l$.
If $yu_3\in E(H)$, then $C'_4$ is an odd hole of length at least $2l+3$, which is not possible. So $C_4\triangle C'_4$ is an odd hole, implying that $C'_4$ is an even hole.

Since $x\neq u_3$, the graph $u_1L_1xP'yP_2u_4$ is a chordal path of $C_2$. Moreover, since $C'_4$ is an even hole, $|Q_2|=1$ by (\ref{eq1}) and Lemma \ref{cz2.3}. Hence, $|P_1|=|L_1|=l$ by (\ref{eq1}) again. 
If $y, u_3$ are non-adjacent, $xP'yP_2u_4Q_2$ is a chordal path of $C_1$, so $xu_1\in E(H)$, implying that $u_4Q_2u_1L_1xP'y$ is a chordal path of $C_3$. Then $yu_4\in E(H)$. 
Since $C_4$ and $C_4\triangle C'_4$ are odd holes, $|P'|=3$, so $|C'_4|=6$, which is not possible. So $yu_3\in E(H)$.
Since $C_4\triangle C'_4$ is an odd hole, $|P'|\geq l+1$ by (\ref{eq1}).
When $x, u_3$ are non-adjacent, since $C'_4$ is an even hole, $C'_4\triangle C_3$ is an odd hole of length at least $2l+3$, which is not possible.
So $xu_3\in E(H)$, implying $|P'|=2l-1$ as $C_4\triangle C'_4$ is an odd hole.
Hence Claim 3.1.4 holds.\qed 

\medskip

\noindent \textbf{Claim 3.1.5} Assume that $P'$ has the structure as stated in Claim 3.1.4. Then no vertex in $V(G) \backslash V(H\cup P')$ has two neighbours in $H\cup P'$.

\pf Assume that a vertex $x''\in V(G) \backslash V(H\cup P')$ has two neighbours $x_1, x_2$ in $H\cup P'$. Since no vertex has two neighbours in an odd hole, 
it follows from Lemma \ref{cz2.7} (3) that $x''$ has exactly two neighbours in $H\cup P'$ with $x_1\in V(H) \backslash \{x, y, u_3\}$ and $x_2\in V(P)$. 
If $x_1\in V(H)\backslash \{P^*_2\cup \{x, u_3\}\}$, then $x_2=v_1$ by Claims 3.1.1-3.1.4.
Now, $G[H\backslash \{P^*_2\}\cup \{x_2, x''\}]$ induces a hole with length $\leq 2l$, a contradiction.
If $x_1\in {P^*_2}\backslash \{y\}$, then by Claims 3.1.1-3.1.4, $x_2=v_2$.
Now, $C_3\cup \{x_2, x''\}$ induces a hole with length $\leq 2l$, a contradiction. \qed

\medskip


By Claims 3.1.1-3.1.3, it suffices to show that $x\notin V(L^*_1\cup L^*_2\cup Q^*_1\cup Q^*_2)$. 
Suppose that it is false.
By symmetry we may assume that $x\in V(L^*_1)$. By Claim 3.1.4, we have that 
$$xu_3\in E(L_1), e=yu_3, |P'|=2l-1, |P_1|=|L_1|=l, |Q_1|=1.$$
Since no $4$-vertex-critical graph has a $P_3$-cut by Lemma \ref{cs}, it suffices to show that $\{x, y, u_3\}$ is a $P_3$-cut of $G$. Assume not. Let $R$ be a shortest induced path in $G-\{x, y, u_3\}$ linking $P$ and $H \backslash \{x, y, u_3\}$.
By 3.1.5, $|R|\geq 3$ and no vertex in $V(H\cup P') \backslash \{x, y, u_3\}$ has a neighbour in $R^*$. Let $s$ and $t$ be the ends of $R$ with $s\in V(P)$.

We claim that $t\notin V(L_1\cup P_2) \backslash \{x, y, u_3\}$. 
Assume to the contrary that $t\in V(L_1) \backslash \{x, u_3\}$ by symmetry.
Let $R_1$ be an induced $(y, t)$-path in $G[V(P'\cup R) \backslash \{x\}]$.
When $u_3$ has no neighbour in $R^*_1$, let $R_2:=R_1$ and $C:=yR_2tL_1u_3y$.
When $u_3$ has a neighbour in $R^*_1$, let $t'\in V(R^*_1)$ be a neighbour of $u_3$ closest to $t$, $R_2:=u_3t'R_1t$ and $C:=u_3R_2tL_1u_3$.
Note that $C_4\triangle C$ is a hole, but $C$ may not be a hole.
Since $C\triangle C_1\triangle C_2$ is an odd hole with length at least $2l+3$ (since $|tL_1u_3|\leq l-1$)
when $C$ is an odd cycle, it suffices to show that $|C|$ is odd. When $x$ has a neighbour in $R^*_2$, since $|L_1|=l$, $|R_2|\geq 2l$. Now $|C_4\triangle C|\geq 3l$,  the subgraph $C_4\triangle C$ is an even hole, which implies  that $C$ is an odd cycle.
So we may assume that $x$ has no neighbour in $R^*_2$.
When $u_3$ is an end of $R_2$, since $R_2$ is a chordal path of $C_1$,
it follows from Lemma \ref{cz2.3} and (\ref{eq1}) that $C$ is an odd hole. 
When $y$ is an end of $R_2$, since neither $x$ nor $u_3$ has a neighbour in $R^*_2$, we have $s=v_1$, by Claims 3.1.1-3.1.4, so $|R_2|>2l$ since $|yP'v_1|=2l-2$ and $|R|\geq 3$.
Then $C_4\triangle C$ is an even hole, so $C$ is an odd cycle.
Hence, the claim holds.

Then $t\in V(P_1\cup L_2) \backslash \{u_1\}$, by symmetry we may therefore assume that $t\in V(P_1) \backslash \{u_1\}$. Let $R_1$ be the induced $(y, t)$-path in $G[V(P'\cup R) \backslash \{x\}]$.
By Claims 3.1.1-3.1.2, either $s=v_1$ and $y$ has no neighbour in $R$ or some vertex in $\{x, u_3\}$ has a neighbour in $R^*_1$.
No matter which case happens, we have $|R_1|\geq 2l$ (if the first case happens, due to $|yP'v_1|=2l-2$ and $|R|\geq 3$, $|R_1|\geq 2l$; else some vertex in $\{x, u_3\}$ has a neighbour in $R^*_1$, and then by $g(G)=2l+1$, $|R_1|\geq 2l$). 
That is, $tR_1yP_2u_4$ is a chordal path of $C_2$ with length at least $3l-1$, which is a contradiction to Lemma \ref{cz2.3} as $t, u_4$ are non-adjacent.
Hence, $\{x, y, u_3\}$ is a $P_3$-cut of $G$.\qed

\begin{theorem} \label{odd theorem}
Let $l\geq 4$ be an integer. For each graph $G$ in ${\cal G}_l$, if $G$ is $4$-vertex-critical, then $G$ has no odd $K_4$-subdivisions.
\end{theorem}
\pf Suppose not. Let $H$ be a subgraph of $G$ that is isomorphic to an odd $K_4$-subdivision.
By Lemma \ref{cz2.7} (2), $H$ is an induced subgraph of $G$. Let $H=(u_1,u_2,u_3,u_4,P_1,P_2,Q_1,Q_2,L_1,L_2)$ (see Figure \ref{fig:$H$}) and
$C_1:=P_1\cup Q_1\cup L_1$, $C_2:=P_1\cup Q_2\cup L_2$, $C_3=P_2\cup Q_1\cup L_2$ and $C_4:=C_1\triangle C_2\triangle C_3$. Since $H$ is an odd $K_4$-subdivision, $C_1$, $C_2$, $C_3$ and $C_4$ are odd holes. By Lemma \ref{cz2.7} (1), 
\begin{equation}\label{eq11}\tag{3.4}
    |P_1|=|P_2|\leq l, |Q_1|=|Q_2|\leq l, |L_1|=|L_2|\leq l. 
\end{equation}
Without loss of generality we may assume that $P_1$, $P_2$ are longest arrises in $H$.

By Lemma \ref{odd lem} and Lemma \ref{cz2.7} (3), any minimal path $P:=(v_1, v_2)$-path linking ${P^*_2}$ and $H \backslash V({P^*_2})$ must have $N_{H \backslash {P^*_2}}(v_1)=\{x\}$, $N_{{P^*_2}}(v_2)=\{y\}$, where $x\in\{u_3, u_4\}$ and $xy\in\{e, f\}$.
Let $P'=xv_1Pv_2y$, so $H\cup P'$ is induced by minimality of $P$.

\medskip 
\noindent \textbf{Claim. 3.2.1} $|P_1|$=$|P_2|=3$.

\pf Suppose not. Since $P_1$, $P_2$ are longest arrises in $H$, $|P_1|$=$|P_2|\geq 4$.
By Lemma \ref{odd lem}, there is a minimal vertex cut $X$ of $G$ with $\{u_3, u_4\}\subseteq X\subseteq N_G[\{u_3, u_4\}]$ and $\{u_3, u_4\}=X\cap V(H)$.
Let $G_1$ be the induced subgraph of $G$ whose vertex set consists of $X$ and the vertex set of the component of $G-X$ containing $P^*_2$.
If all induced $(u_3, u_4)$-paths in $G_1$ have length $|P_2|$, by Lemma \ref{cz2.4}, $G$ has a degree-$2$ vertex or a $K_2$-cut, which is not possible.
Hence, it suffices to show that all induced $(u_3, u_4)$-paths in $G_1$ have length $|P_2|$.

Let $Q$ be an arbitrary induced $(u_3, u_4)$-paths in $G_1$. 
When $|L_1|\geq 2$, since $QQ_2$ is a chordal path of $C_1$ by Lemma \ref{cz2.7} (3) and the definition of $G_1$, we have $|QQ_2|=|Q_1P_1|$ by Lemma \ref{cz2.3}, so $|Q|=|P_1|$ by (\ref{eq11}).
Hence, by (\ref{eq11}) we may assume that $|L_1|=1$ and  $|Q_1|=|P_1|=l$.
Since $Q_1L_2$ is an induced $(u_3, u_4)$-path  of length $l+1$, 
either $|Q|=|P_2|=l$ or $|Q|\geq l+1$ and $|Q|$ has the same parity as $l+1$.
Assume that the latter case happens, then $G[L_1\cup P_1\cup L_2\cup Q]$ is an odd hole of length at least $2l+3$, which is not possible.\qed

\medskip

By Claim 3.2.1, $G\in {\cal G}_4$, so $g(G)=9$.
Thus $H$ is a regular odd $K_4$-subdivision  and each arris of $H$ has length 3.
And by Lemma \ref{odd lem}, without loss of generality, we may assume $x=u_3$, then $xy=e$.
Since $4$-vertex-critical graph $G$ has no $P_3$-cut (in particular $\{x,y,v_1\}$ is not a $P_3$-cut), there exists a direct connection induced $(w,t)$-path $R'$ linking $P\backslash\{v_1\}$ and $H\backslash \{x, y\}$ in $G-\{x, y, v_1\}$ with $N_{H\backslash \{x, y\}}(w) \ne \emptyset$.
It is clear that $|N_H(w)|=|N_{H\backslash \{x, y\}}(w)|=1$.
Let $R$ be the path induced by $R'\cup (P\backslash \{v_1\})\cup \{y\}\cup N_{H}(w)$, so one end of $R$ is $y$ and the other is in $N_{H}(w)$.
Let the neighbour of $y$ in $P^*_2$ is $y'$.
Note that $u_3$ has a neighbour in $R$ by Lemma \ref{odd lem} as $u_3 \not\in N_H(w)$.

We claim that $wy'\notin E(G)$.
Otherwise, there exists a path $Q$ induced by $P'\cup R\backslash \{y\}$ linking $u_3$ and $y'$.
Set $C'=u_3L_1u_1P_1u_2L_2u_4y'Qu_3$, $C''=u_3L_1u_1Q_2u_4y'Qu_3$.
Then $|C'|=10+|Q|$ and $|C''|=7+|Q|$ imply $|Q|=2$, which contradicts Lemma \ref{cz2.7} (3).

We claim that $wu_1, wu_2\notin E(G)$.
Otherwise, we may assume $wu_1\in E(G)$.
Set $C'=u_1Q_2u_4y'yRu_1$, $C''=u_1P_1u_2L_2u_4yy'Ru_1$.
Then $|C'|=5+|R|$ and $|C''|=8+|R|$ imply $|R|=4$.
Now, $G[\{u_3, y\}\cup R]$ contains a cycle with length less than $9$, a contradiction.

We claim that $N_{ P^*_1}(w)=\emptyset$.
Otherwise, let $u'_1$ be the neighbour of $u_1$ in $P_1$.
Without loss of generality, we may assume that $wu'_1\in E(G)$.
Set $C'=u'_1u_1Q_2u_4y'yRu'_1$, $C''=u'_1P_1u_2L_2u_4y'yRu'_1$.
Then $|C'|=6+|R|$ and $|C''|=7+|R|$, implying that $|R|=3$.
Now, $G[\{u_3, y\}\cup R]$ contains a cycle with length less than $9$, a contradiction.

We claim that $wu_4\notin E(G)$.
Otherwise, let $u'_3$ be the neighbour of $u_3$ in $R$ closest to $w$.
Set $C'=u'_3Ru_4Q_2u_1L_1u_3u'_3$, $C''=u'_3Ru_4Q_2u_1P_1u_2Q_1u_3u'_3$.
Then $|C'|=7+|u'_3Ru_4|$ and $|C''|=10+|u'_3Ru_4|$ imply that $|u'_3Ru_4|=2$.
Since $yu'_3\notin E(G)$, now, $yu_3u'_3Ru_4y'y$ is a $6$-hole, a contradiction.

We claim that $N_{Q^*_1\cup L^*_1}(w)=\emptyset $.
Otherwise, Let $N_{H}(w)=\{z\}$ and without loss of generality
we may assume that $z\in Q^*_1$.
Let $h := |zQ_1u_2|$.
Note that $h\in \{1,2\}$.
Set $C'=zQ_1u_2L_2u_4y'yRz$, $C''=zQ_1u_2P_1u_1Q_2u_4y'yRz$.
Then $|C'|=h+5+|R|$ and $|C''|=h+8+|R|$.
So $|R| = 4 - h \le 3$.
Now,
$G[\{u_3, y\}\cup R]$ contains a cycle with length less than $9$, a contradiction.

Hence, we have $N_{Q^*_2\cup L^*_2}(w) \ne \emptyset $.
Let $N_{H}(w)=\{z\}$
and without loss of generality, we may assume that $z\in Q^*_2$.
If $zu_1\in E(G)$, set $C'=zu_1P_1u_2L_2u_4y'yRz$, $C''=zQ_2u_4y'yRz$.
Then $|C'|=9+|R|$ and $|C''|=4+|R|$, implying $|R|=5$.
Now, $G[\{u_3, y\}\cup R]$ contains a cycle with length less than $9$, a contradiction.
So $zu_4\in E(G)$, let $u'_3$ be the neighbour of $u_3$ in $R$ closest to $w$.
Then $u'_3y\notin E(G)$.
Set $C'=u'_3Rzu_4y'yu_3u'_3$ and 
$C''=u'_3RzQ_2u_1L_1u_3u'_3$.
Then $|C'|=|u'_3Rz|+5$ and $|C''|=|u'_3Rz|+6$, implying $|u'_3Rz|=4$.
Now, $u'_3RzQ_2u_1P_1u_2Q_1u_3u'_3$ is a $13$-hole, a contradiction.
\qed

\begin{theorem} \label{theorem2}
Let $G\in {\cal G}_4$. Assume $G$ has no 2-edge-cut or $K_2$-cut.
Then one of the following holds.
\begin{enumerate}[1)]
    \item $G$ has an odd $K_4$-subdivision.
    \item $G$ contains a balanced $K_4$-subdivision of type $(1, 2)$.
    \item $G$ has a $P_3$-cut.
    \item $G$ has a degree-2 vertex.
\end{enumerate}
\end{theorem}
\pf Assume that neither 1) nor 2) is true.
Set $C=v_1v_2\cdots v_9v_1$. 
Since $G$ has no $2$-edge cut or $K_2$-cut, there exists a jump over $C$.
By Corollary \ref{coro-jump}, $C$ either has a short jump or a local jump over $C$ across one vertex.
If there exists a short jump over $C$, let $P$ be a short jump over $C$ with $|P|$ as small as possible.
By symmetry we may assume that the ends of $P$ are $v_2, v_k$ with $k\leq 6$ and $P$ is across $v_2v_3 \ldots v_k$ (as $g(G) = 9$).

By  Lemma \ref{c234}, it is clear that any two short jumps are not crossing.
Then, if there exists a second short jump $P^*$, then either one of the ends of $P^*$ is in $\{v_2, v_k\}$ and the other end is in $\{v_{k+1}, \cdots, v_9\}$, or both ends of $P^*$ are in $\{v_{k+1}, \cdots, v_9\}$,
then either $|P^*|<|P|$ or  by Lemma \ref{c235} 1), $G$ has an odd $K_4$-subdivision, a contradiction.
So all short jumps over $C$ have ends in $\{v_2, v_3, \cdots, v_k\}$ and are across a subpath of $C(v_2,v_k)$.
Since $|P|$ is minimum, we have additionally that
\begin{equation}\label{eq4}\tag{3.5}
\mbox{no jump hole over $C$ contains a vertex in $V(C)\backslash\{v_2, v_3, \cdots, v_k\}$.}
\end{equation}

By Lemma \ref{c235} 2), 
as long as there exist two local jumps over $C$ across one vertex that are uncrossing, then there exists a jump hole, and thus a short jump.
Otherwise, all local jumps over $C$ across one vertex are crossing or there exists only one local jump over $C$ across one vertex.
For each integer $1\leq i\leq 2$, let $P_i$ be a local $(s_i, t_i)$-jump over $C$ across one vertex and $Q_i$ be the $(s_i, t_i)$-path on $C$ of length $2$.
When $P_1$ and $P_2$ are uncrossing, by (\ref{eq4}) and Lemma \ref{c235} 2),
$|V(Q_1\cup Q_2)\backslash\{v_2, v_3, \cdots, v_k\}|\leq 2$,
 then at least one vertex in $\{s_i, t_i\}$ for each $i \in \{1,2\}$ is in $\{v_2, v_3, \cdots, v_k\}$.
When $P_1$ and $P_2$ are crossing, $|V(Q_1\cup Q_2)|=4$.
If there is only one local jump over $C$ across one vertex denoted by $P_1$, then $V(Q_1)=3$ where $P_1$ across $Q_1$ over $C$.
\begin{figure}[h!] 
\centering
\includegraphics[width=0.7\linewidth]{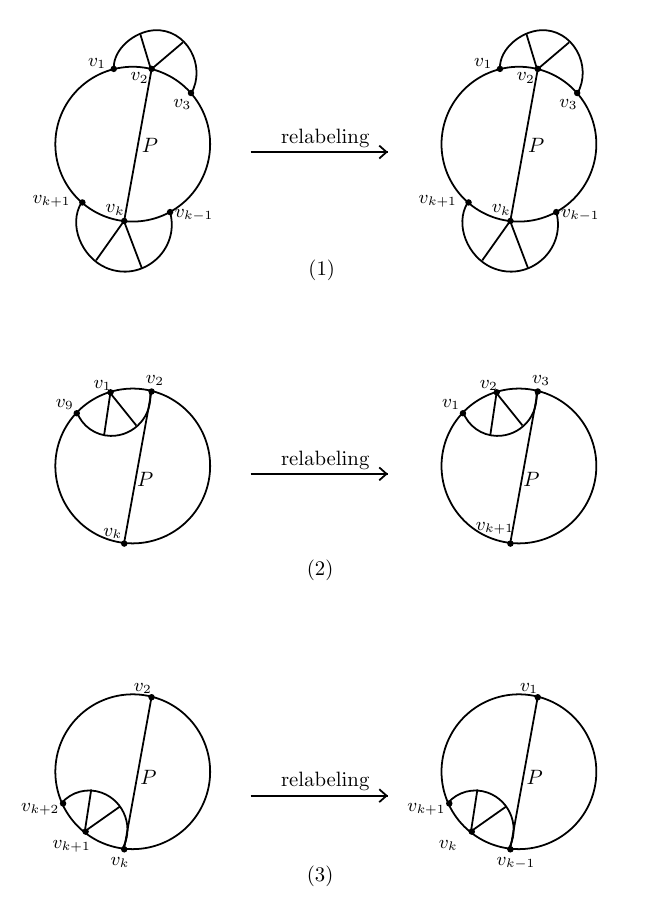}
\caption{\label{relabeling} \centering Relabeling $V(C)$ \\
     (1) one local jump across $v_2$ and one local jump across $v_k$, no relabeling \\
     (2) a local jump across $v_1$, relabel $v_i$ to be $v_{i+1}$  \\
     (3) a local jump across $v_{k+1}$, relabel $v_i$ to be $v_{i-1}$
}
\end{figure}

We relabel the indices of $v_1, \ldots, v_9$ for convenience. 
If there is no short jump, we may assume the ends of all local jump across one vertex are in $\{v_1, v_2, \cdots, v_{k+1}\}$.
Now assume $P$ exists.
If there exists a local jump across $v_1$, then we relabel $v_i$ to be $v_{i+1}$ for each $i\in [9]$ (we write $v_9$ instead of $v_0$ for convenience), see Figure \ref{relabeling} (2) for illustration.
And if there exists a local jump across $v_{k+1}$, then we relabel $v_i$ to be $v_{i-1}$ for each $i\in [9]$ (we write $v_9$ instead of $v_0$ for convenience), see Figure \ref{relabeling} (3) for illustration.
For the rest of cases, we do not relabel.
We would like to point out that if there exist two local jump with one across $v_2$ and the other across $v_k$, then we do not relabel, see Figure \ref{relabeling} (1) for illustration.

In all the cases, after possible relabeling, we may assume that 
\begin{equation}\label{eq5}\tag{3.6}
\begin{split}
    &\mbox{all local jumps over $C$ across one vertex and all short jumps have ends in $\{v_1, v_2, \cdots, v_k, v_{k+1}\}$ } \\
    &\mbox{with $4\leq k\leq 6$ and all short jumps are across a subpath of $v_1v_2 \cdots v_k$ of length at most $k-2$.}
\end{split}
\end{equation}

For any integer $1\leq i\leq k+1$, let $X_i$ be the set of vertices adjacent to $v_i$ that are in a local jump over $C$ across one vertex with one end $v_i$ or a short jump over $C$ with one end $v_i$. Set $X=X_1\cup X_2\cup \cdots\cup X_{k+1}$.

Since no vertex in $V(G)\backslash V(C)$ has two neighbours in $V(C)$, $v_8$ has no neighbours in $X$. Assume that $v_8$ has degree at least 3, for otherwise 4) holds. There is a connected induced subgraph $D$ such that $v_8$ has a neighbour in $V(D)$ and $V(D)\cap (V(C)\cup X)=\emptyset$, and we choose $D$ to be maximal with these properties.
Let $N=\{w|w\in V(C)\cup X, N_D(w)\neq \emptyset\}$.
Evidently, $v_8\in N$.

\medskip

It suffices to show that $N\subseteq \{v_7, v_8, v_9\}.$ 
Because otherwise $\{v_7, v_8, v_9\}$ is a $P_3$-cut, and thus 3) holds. 
This will conclude the proof.

Suppose that $N \not\subseteq \{v_7, v_8, v_9\}$. For $1\leq i\leq k+1\leq 7$, let $W_i=X_i\cup \{v_i\}$.
First assume that $N\cap W_i\neq \emptyset$ for some $i\in [6]$.
Let $Q$ be a shortest $(v_8, v_i)$-path with interior in $V(D)\cup (W_i\backslash\{v_i\})$.
Exactly one of $\{v_7, v_9\}$ has a neighbour in $Q^*$, otherwise either there is a $(v_7, v_9)$-short jump or $Q$ is a short jump,  which contradicts to (\ref{eq5}).

We claim that $Q$ is not a $(v_8, v_1)$-path.
Suppose not. If $N_{Q^*}(v_7)=\emptyset$, then $N_{Q^*}(v_9)\neq\emptyset$ and $Q$ is a local jump across $v_9$,  which contradicts to (\ref{eq5}).
Otherwise, there is a $(v_1, v_7)$-short jump, a contradiction to (\ref{eq5}).

So it is clear that $Q$ is a $(v_8, v_i)$-path for $i\in \{2,3,\cdots,6\}$.
If $N_{Q^*}(v_7)=\emptyset$ and $N_{Q^*}(v_9)\neq\emptyset$, then there is a $(v_i, v_9)$-short jump,  which contradicts (\ref{eq5}).
Thus we have
\begin{equation}\label{eq6}\tag{3.7}
 \mbox{ $N_{Q^*}(v_7)\neq\emptyset$ and $N_{Q^*}(v_9)=\emptyset$ for all $i\in \{2,3,\cdots,6\}$. }
\end{equation}

We claim that $Q$ is not a $(v_8, v_2)$-path.
Otherwise, by (\ref{eq6}), there is a $(v_2, v_7)$-short jump, which contradicts (\ref{eq5}).

\medskip
\noindent \textbf{Claim.} $Q$ is not a $(v_8, v_i)$-path for $i\in \{3, 4, 5\}$.

\pf By (\ref{eq6}), there is a $(v_i, v_7)$-short jump, which implies that $C$ is relabelled. 
Then there is a local jump across only $v_2$, say $R$.
So $|R|$ is odd with length $\geq 7$ by Lemma \ref{c236}.

Let $Q':=(v_i,v_7)$-short jump, where $Q'\backslash \{v_7\}\subset Q$, then $|Q'|=9-(7-i)=2+i$. Let $N_R(v_3)=\{x\},$ $N_{Q'}(v_i)=\{y\}$.

If there exists a vertex in $V(R^*)\cap V(Q'^*) \setminus \{x,y\}$ or the end of the edge between $R^*$ and ${Q'}^*$ is not $x$ or $y$, then there is a $(v_7, v_1)$-short jump or $(v_7, v_2)$-short jump, which is a contradiction as $|P|$ is minimum. Then we have
\begin{enumerate}
\setlength{\itemsep}{0pt}
    \item[(i)] either $V(R^*)\cap V(Q'^*)=\emptyset$ and $R^*$ is anticomplete to $Q'^*$ 
    \item[(ii)] or $N_{Q'^*}(x)\neq\emptyset$
    \item[(iii)]  or $N_{R^*}(y)\neq\emptyset$
    \item[(iv)]  or $V(R^*)\cap V(Q'^*)\cap \{x, y\}\neq\emptyset$
\end{enumerate}
The first case implies an odd hole $v_1Rv_3Cv_iQ'v_7v_8v_9v_1$ of length $\geq 11$, a contradiction.
The second case is impossible since $|Q|$ is shortest and $g(G)=9$.
If the third case happens, let $z\in N_R(y)$ be nearest to $v_1$. 
When $i=3$, since $v_1Rzyv_3$ is a jump across $v_2$, $|v_1Rz|$ is odd and has size at least $7$. Then $v_7Q'yzRv_1v_9v_8v_7$ is an odd hole with length $8+|v_1Rz| \ge 15$, a contradiction. 
When $i=\{4, 5\}$, then $v_1Rzyv_i$ is either a short jump (which is a contradiction) or a local jump has the same parity as $Q'$.
Now, $v_1v_9v_8v_7Q'yRzv_1$ is an odd hole with length $\geq 11$, a contradiction.
So $V(R^*)\cap V(Q'^*)\cap \{x, y\}\neq\emptyset$.
Thus either $x \in V(R^*)\cap V(Q'^*)$ or $y \in V(R^*)\cap V(Q'^*)$.
In both cases, one can see that $x=y$ and $i=3$.
Similarly, there is an odd hole with length $8+|v_1Rz| \ge 15$, a contradiction.
\qed

\medskip

We claim that $Q$ is not a $(v_8, v_6)$-path.
Otherwise, by (\ref{eq6}), there is a $(v_8, v_6)$-local jump across $v_7$, which contradicts to (\ref{eq5}).

Therefore, such $Q$ does not exist, so there is no path from $v_8$
to $\{v_1, v_2,\cdots, v_6\}$ with interior in $V(D)\cup (W_i\backslash\{v_i\})$.

\medskip

If $X_7\cap N\neq \emptyset$, then there exists $x_7\in X_7\cap N$. 
Since $x_7$ is in a local jump or a short jump, there is a path from $v_8$
to $\{v_1, v_2,\cdots, v_6\}$, a contradiction. 
So $X_7\cap N=\emptyset$.
This completes the proof.
\qed


\medskip 

Now we are ready to prove Theorem \ref{main theorem}.

\pf (Theorem \ref{main theorem})
Suppose Theorem \ref{main theorem} is not true. 
Let $G\in {\cal G}_4$ be a minimal counterexample to Theorem \ref{main theorem}.
Then $G$ is $4$-vertex-critical. So $G$ has no degree-2 vertex or 2-edge cut or $K_2$-cut.
By Theorem \ref{odd theorem} and Lemma \ref{c231}, $G$ contains neither odd $K_4$-subdivision nor balanced $K_4$-subdivision of type $(1, 2)$. Hence, $G$ has a $P_3$-cut by Theorem \ref{theorem2}, which is a contradiction to Lemma \ref{cs}.
\qed

\end{document}